\documentclass{amsart}

\usepackage{amsmath,amsthm,amsopn,amstext,amscd,amsfonts,amssymb,mathrsfs,mathtools}
\usepackage{dsfont}
\usepackage{comment}
\usepackage{xcolor}
\usepackage{float}
\usepackage[active]{srcltx}
\usepackage{graphicx, mathdots}

  \usepackage{mathtools}

\newtheorem{theorem}{\sc Theorem}[section]
\newtheorem{conj}{\sc Conjecture}[section]

\usepackage[active]{srcltx}

\usepackage{graphicx, epsfig, subfig}
\usepackage{lscape}
\usepackage{rotating}
\usepackage{caption}
\usepackage{multicol}
\usepackage{colortbl}
\usepackage{nicematrix}

\newmuskip\pFqskip
\mathchardef\pFcomma=\mathcode`, 

\makeatletter
\def\BState{\State\hskip-\ALG@thistlm}
\makeatother

\def\downbar#1{
\setbox10=\hbox{$#1$}
            \dimen10=\ht10 \advance\dimen10 by 2.5pt
            \ifdim \dimen10<15pt 
               \advance\dimen10 by -0.5pt
               \dimen11=\dimen10
               \advance\dimen10 by 2.5pt
               \lower \dimen11
            \else \lower \ht10 \fi
            \hbox {\hskip 1.5pt \vrule height \dimen10 depth \dp10}}
\def\upbar#1{
\setbox10=\hbox{$#1$}
            \dimen10=\ht10 \advance\dimen10 by \dp10 \advance\dimen10 by 2.5pt
            \ifdim \dimen10<15pt 
                \advance\dimen10 by 2pt \fi
            \raise 2.5pt \hbox {\hskip -1.5pt \vrule height \dimen10}}

\begin{document}
\title[Proof of two conjectures on Askey-Wilson polynomials]{Proof of two conjectures on Askey-Wilson polynomials}
\author{K. Castillo}
\address{CMUC, Department of Mathematics, University of Coimbra, 3001-501 Coimbra, Portugal}
\email{ kenier@mat.uc.pt}
\author{D. Mbouna}
\address{CMUC, Department of Mathematics, University of Coimbra, 3001-501 Coimbra, Portugal}
\email{dmbouna@mat.uc.pt}

\date{\today}
\begin{abstract}
We give positive answer to two conjectures posed by M. E. H Ismail in his monograph [Classical and quantum orthogonal polynomials in one variable, Cambridge University Press, 2005].
\end{abstract}
\maketitle
\section{Introduction and main result}
The Askey-Wilson divided difference operator is defined by
\begin{align}
(\mathcal{D}_q  f)(x)=\frac{\breve{f}\big(q^{1/2} z\big)
-\breve{f}\big(q^{-1/2} z\big)}{\breve{e}\big(q^{1/2}z\big)-\breve{e}\big(q^{-1/2} z\big)},\quad
z=e^{i\theta}, \label{0.3}
\end{align}
where $\breve{f}(z)=f\big((z+1/z)/2\big)=f(\cos \theta)$ for each polynomial $f$ and $e(x)=x$.
Here $0<q<1$ and $\theta$ is not necessarily a real number (see \cite[p. 300]{I05}). Hereafter, we denote $x(s)=(q^{s}+q^{-s})/2$ with $0<q<1$. Taking $e^{i\theta}=q^s$ in \eqref{0.3}, $\mathcal{D}_q$ reads
\begin{align*}
\mathcal{D}_q f(x(s))= \frac{f\big(x(s+\frac{1}{2})\big)-f\big(x(s-\frac{1}{2})\big)}{x(s+\frac{1}{2})-x(s-\frac{1}{2})}.
\end{align*}
Set $\mathcal{D}^0_q\, f=f$ and $\mathcal{D}^1_q=\mathcal{D}_q$, and define $\mathcal{D}^k_q=\mathcal{D}_q(\mathcal{D}^{k-1}_q)$ for each $k=1,2, \dots$. The following two conjectures, which generalize the Sonin-Hahn problem, were posed by M. E. H Ismail in his monograph on Orthogonal Polynomials and Special Functions published in 2005 (see \cite[Conjecture 24.7.10 and Conjecture 24.7.11]{I05})  and revised in 2009:

\begin{conj}\label{C1}
If $(p_n)_{n\geq 0}$ and $(\mathcal{D}_q p_{n})_{n\geq 0}$, or the latter with a limiting case of $\mathcal{D}_q$, are two sequences of orthogonal polynomials, then $(p_n)_{n\geq 0}$ are multiples of the Askey-Wilson polynomials, or special or limiting cases of them.
\end{conj}

\begin{conj}\label{C2}
If $(p_n)_{n\geq 0}$ and $(\mathcal{D}^k_q p_{n+k})_{n\geq 0}$, or the latter with a limiting case of $\mathcal{D}_q$, are two sequences of orthogonal polynomials for some $k$, $k=1,2,\dots$, then $(p_n)_{n\geq 0}$ are multiples of the Askey-Wilson polynomials, or special or limiting cases of them.
\end{conj}

Define the average operator $\mathcal{S}_q$ by
\begin{align*}
\mathcal{S}_q f(x(s))=\frac{f\big(x(s+\frac{1}{2})\big)+f\big(x(s-\frac{1}{2})\big)}{2}.
\end{align*}
 for every polynomial $f$. In 2003, Ismail proved the following result (see \cite[Theorem 20.1.3]{I05}):
\begin{theorem}\label{T}
A second order operator equation of the form
\begin{align}\label{ismail}
f(x)\mathcal{D}^2_q\, y+g(x) \mathcal{S}_q \mathcal{D}_q \, y+h(x)\, y=\lambda_n\, y
\end{align}
has a polynomial solution $y_n(x)$ of exact degree $n$ for each $n=0,1,\dots$, if and only if $y_n(x)$ is a multiple of the Askey-Wilson polynomials, or special or limiting cases of them. In all these cases $f$, $g$, $h$, and $\lambda_n$ reduce to
\begin{align*}
f(x)&=-q^{-1/2}(2(1+\sigma_4)x^2-(\sigma_1+\sigma_3)x-1+\sigma_2-\sigma_4),\\[7pt]
g(x)&=\frac{2}{1-q} (2(\sigma_4-1)x+\sigma_1-\sigma_3), \quad h(x)=0,\\[7pt]
\lambda_n&=\frac{4 q(1-q^{-n})(1-\sigma_4 q^{n-1})}{(1-q)^2}, 
\end{align*}
or a special or limiting case of it, $\sigma_j$ being the jth elementary symmetric function of the Askey-Wilson parameters. 
\end{theorem}

Virtually the above conjectures are summed up in one if we are able to prove Conjecture \ref{C2}. To do this, we prove that the sequences of polynomials appearing in Conjecture \ref{C2} satisfy, for each $k$, a second order operator equation of the form \eqref{ismail}. The important point to note here is that this argument would not lead to a satisfactory conclusion if we were not looking for the whole space of ``Askey-Wilson polynomials, or special or limiting cases of them".

\begin{theorem}\label{L}
If $(p_n)_{n\geq 0}$ and $(\mathcal{D}^k_q p_{n+k})_{n\geq 0}$, or the latter with a limiting case of $\mathcal{D}_q$, are two sequences of orthogonal polynomials for some $k$, $k=1,2,\dots$, then, for each $k$, $(\mathcal{D}^{k-1}_q p_{n+k-1})_{n\geq 0}$ are multiples of the Askey-Wilson polynomials, or special or limiting cases of them.
 \end{theorem}
 
 Fix $k$, $k=1,2,\dots$. It is easily seen that $(p_n)_{n\geq 0}$ is a sequence of orthogonal polynomials satisfying 
 \begin{align}\label{estructura}
 \pi(x) \mathcal{D}^k_q p_{n}(x)=\sum_{j=-m}^m c_{n,j}\, p_{n+j}(x),\quad c_{n,-m}\not=0,
 \end{align}
for a polynomial $\pi$ which does not depend on $n$, if and only if $(p_n)_{n\geq 0}$ and $(\mathcal{D}^k_q p_{n+k})_{n\geq 0}$ are sequences of orthogonal polynomials. Now, we can apply Theorem \ref{L} to conclude that for each $k$, $k=1,2,\dots$, $(\mathcal{D}^{k-1}_q p_{n+k-1})_{n\geq 0}$ are multiples of the Askey-Wilson polynomials, or special or limiting cases of them. In particular, taking $k=2$ and $m=2$ in \eqref{estructura}, we have the main result proved in \cite{KJ19}: $(p_n)_{n\geq 0}$ are multiples of the Askey-Wilson polynomials, or special or limiting cases of them. Neither in this work nor in \cite{KJ19} was possible to exclude the `limiting cases' in the last statement. If so, we would have positive answer to a particular case of another conjecture posed by Ismail (see \cite[Conjecture 24.7.9]{I05}).
\section{Proof of Theorem \ref{L}}
The following properties are well known:
\begin{align*}
\mathcal{D}_q (fg)&= (\mathcal{D}_q f)(\mathcal{S}_q g)+(\mathcal{S}_q f)(\mathcal{D}_q g),\\[7pt]
\mathcal{S}_q ( fg)&=\mathrm{U}_2 (\mathcal{D}_q f) (\mathcal{D}_q g)  +(\mathcal{S}_q f) (\mathcal{S}_q g),\\[7pt]
\mathcal{S}^2_q f&=\alpha  \mathrm{U}_2 \mathcal{D}^2_q f+\mathrm{U}_1 \mathcal{S}_q \mathcal{D}_q f+f,\\[7pt]
\mathcal{D}_q \mathcal{S}_q f&=\alpha \mathcal{S}_q \mathcal{D}_q f+\mathrm{U}_1\mathcal{D}^2_q f,
\end{align*}
for polynomials $f$ and $g$, where $\mathrm{U}_1(x)=(\alpha^2-1)x$, $\mathrm{U}_2(x)=(\alpha^2-1)(x^2-1)$, and  $2 \alpha=q^{1/2}+q^{-1/2}$.
We leave it to the reader to verify by induction that
 \begin{align*}
 \mathcal{D}^k_q(X f(x))=\gamma_k \mathcal{S}_q \mathcal{D}^{k-1}_q f(x)+\frac{q^{k/2}+q^{-k/2}}{2} X \mathcal{D}^k_q f(x).
 \end{align*}
where 
 \begin{align*}
 \gamma_k=\frac{q^{k/2}-q^{-k/2}}{q^{1/2}-q^{-1/2}}.
 \end{align*}
Set $P_n^{[k]}=\gamma_n!/\gamma_{n+k}!\,\mathcal{D}^k_q P_{n+k}$, and so $P_n^{[k]}=\gamma_{n+1}^{-1}\mathcal{D}_q P^{[k-1]}_{n+1}$. Since $(P_n)_{n\geq 0}$ and $(P^{[k]}_n)_{n\geq 0}$, for a certain fixed $k$, are sequences of (monic) orthogonal polynomials, any three consecutive elements of these sequences satisfy
\begin{align}
\label{rec1}x P_{n+k}(x)&=P_{n+k+1}(x)+B_{n+k}P_{n+k}(x)+C_{n+k} P_{n+k-1}(x),\\[7pt]
\label{rec2}x P^{[k]}_{n-1}(x)&=P^{[k]}_{n}(x)+B^{[k]}_{n-1}P^{[k]}_{n-1}(x)+C^{[k]}_{n-1} P^{[k]}_{n-2}(x),
\end{align}
with $C_{n+k}\not=0$ and $C^{[k]}_{n-1}\not=0$. We apply $\mathcal{D}^{k}_q$ to \eqref{rec1} to get
\begin{align}
\label{1a}&\gamma_{k}\mathcal{S}_q P^{[k-1]}_{n}(x)+\frac{q^{k/2}+q^{-k/2}}{2} x  \mathcal{D}_q P^{[k-1]}_{n}(x)\\[7pt]
\nonumber &\quad =\frac{\gamma_{n+k}}{\gamma_{n+1}}\mathcal{D}_q P^{[k-1]}_{n+1}(x)+B_{n+k-1}\mathcal{D}_q P^{[k-1]}_{n}(x)+\frac{\gamma_{n}}{\gamma_{n+k-1}}C_{n+k-1}\mathcal{D}_q P^{[k-1]}_{n-1}(x).
\end{align}
From \eqref{rec2} we have
\begin{align}
\label{1b}&\frac{1}{\gamma_n}x\mathcal{D}_q P^{[k-1]}_{n}(x)\\[7pt]
\nonumber &\quad =\frac{1}{\gamma_{n+1}}\mathcal{D}_q P^{[k-1]}_{n+1}(x)+\frac{1}{\gamma_n}B^{[k]}_{n-1} \mathcal{D}_q P^{[k-1]}_{n}(x)+\frac{1}{\gamma_{n-1}}C^{[k]}_{n-1} \mathcal{D}_q P^{[k-1]}_{n-1}(x).
\end{align}
We now apply $\mathcal{S}_q$ to \eqref{1a} and \eqref{1b} and, by combining the resulting equations, we can eliminate $\mathcal{S}_q \mathcal{D}_q  P^{[k-1]}_{n-1}(x)$, and obtain the equation
\begin{align}
&\label{3} D_n(x)  \mathcal{S}_q \mathcal{D}_q  P^{[k-1]}_{n}(x)+E_n \mathrm{U}_2(x) \mathcal{D}^2_qP^{[k-1]}_{n}(x)+\frac{\gamma_k}{\gamma_{n-1}} C_{n-1}^{[k]}P^{[k-1]}_{n}(x)\\[7pt]
\nonumber &\quad =F_n\mathcal{S}_q \mathcal{D}_q  P^{[k-1]}_{n+1}(x),
\end{align}
where 
\begin{align*}
D_n(x)&=\left(\frac{q^{(k+1)/2}+q^{-(k+1)/2}}{2}\frac{1}{\gamma_{n-1}}C_{n-1}^{[k]}-\frac{\alpha}{\gamma_{n+k-1}}C_{n+k-1} \right)x-\frac{1}{\gamma_{n-1}}B_{n+k-1}C^{[k]}_{n-1}\\[7pt]
&\quad +\frac{1}{\gamma_{n+k-1}}B^{[k]}_{n-1}C_{n+k-1},\\[7pt]
E_n&=\frac{\gamma_{k+1}}{\gamma_{n-1}} C_{n-1}^{[k]}-\frac{1}{\gamma_{n+k-1}}C_{n+k-1},\\[7pt]
F_n&=\frac{\gamma_{n+k}}{\gamma_{n+1}\gamma_{n-1}} C^{[k]}_{n-1}-\frac{\gamma_{n}}{\gamma_{n+1}\gamma_{n+k-1}} C_{n+k-1}.
\end{align*}
Similarly, we can eliminate $\mathcal{S}_q \mathcal{D}_q  P^{[k-1]}_{n+1}(x)$, and  shift $n$ to $n+1$ to obtain the equation
\begin{align}
\label{4}&\widetilde{D}_n(x)  \mathcal{S}_q \mathcal{D}_q  P^{[k-1]}_{n+1}(x)-\widetilde{E}_n \mathrm{U}_2(x) \mathcal{D}^2_qP^{[k-1]}_{n+1}(x)-\frac{\gamma_k}{\gamma_{n+2}} P^{[k-1]}_{n+1}(x)\\[7pt]
\nonumber &\quad=F_{n+1}\mathcal{S}_q \mathcal{D}_q  P^{[k-1]}_{n}(x),
\end{align}
where
\begin{align*}
\gamma_{n+2} \widetilde{D}_n(x)&=\frac{q^{k/2}+q^{-k/2}}{2} \frac{\gamma_{n}}{\gamma_{n+1}} x+B_{n+k}-\frac{\gamma_{n+k+1}}{\gamma_{n+1}} B_n^{[k]},\\[7pt]
\widetilde{E}_n&=\frac{\gamma_n \gamma_k}{\gamma_{n+1}\gamma_{n+2}}.
\end{align*}
We now apply $\mathcal{D}_q$ to \eqref{1a} and \eqref{1b} and, by combining the resulting equations, we can eliminate $\mathcal{D}^2_q  P^{[k-1]}_{n-1}(x)$, and obtain the equation
\begin{align}\label{3b}
D_n(x)\mathcal{D}^2_q  P^{[k-1]}_{n}(x)+E_n  \mathcal{S}_q \mathcal{D}_q  P^{[k-1]}_{n}(x)=F_{n} \mathcal{D}^2_q  P^{[k-1]}_{n+1}(x).
\end{align}
Similarly, we can eliminate $\mathcal{D}^2_q  P^{[k-1]}_{n+1}(x)$, and shift $n$ to $n+1$ to obtain the equation
\begin{align}\label{3a}
\widetilde{D}_n(x)\mathcal{D}^2_q  P^{[k-1]}_{n+1}(x)-\widetilde{E}_{n}  \mathcal{S}_q \mathcal{D}_q  P^{[k-1]}_{n+1}(x)=F_{n+1}\mathcal{D}^2_q  P^{[k-1]}_{n}(x).
\end{align}
Note that if $\mathcal{S}_q P^{[k]}_{n}(x(s_1))=0$ we have $P^{[k]}_{n}(x(s_1+1/2))=-P^{[k]}_{n}(x(s_1-1/2))$. Suppose that $\mathcal{D}_q  P^{[k]}_{n}(x(s_1))=0$. So 
$$
P^{[k]}_{n}(x(s_1+1/2))=P^{[k]}_{n}(x(s_1-1/2))=0,
$$
which is impossible. Thus $\mathcal{D}_q  P^{[k]}_{n}(x)$ and $\mathcal{S}_q P^{[k]}_{n}(x)$ have no common zeros. After shifting $n$ to $n-1$ in \eqref{3a},  to obtain a contradiction, suppose that $F_n=0$, i.e. 
$$
\widetilde{D}_{n-1}(x)\mathcal{D}_q  P^{[k]}_{n}(x)=\widetilde{E}_{n-1}  \mathcal{S}_q  P^{[k]}_{n}(x).
$$
Since $\mathcal{D}_q  P^{[k]}_{n}(x)$ and $\mathcal{S}_q P^{[k]}_{n}(x)$ have no common zeros, we have $\widetilde{E}_{n-1}=0$, which is impossible. (We can also conclude that $F_{n+1}\not=0$.) Multiplying \eqref{3a} by $F_n$ and using \eqref{3b} and \eqref{3}, we get
\begin{align}\label{6}
f_n(x)\mathcal{D}^2_q  P^{[k-1]}_{n}(x)+g_n(x) \mathcal{S}_q \mathcal{D}_q  P^{[k-1]}_{n}(x)+\frac{\gamma_k}{\gamma_{n-1}}\widetilde{E}_n C^{[k]}_{n-1}P^{[k-1]}_{n}(x)=0,
\end{align}
where
\begin{align*}
f_n(x)&= E_n \widetilde{E}_n\mathrm{U}_2-D_n(x) \widetilde{D}_n(x)+F_n F_{n+1},\\[7pt]
 g_n(x)&=\widetilde{E}_n D_n(x)-E_n \widetilde{D}_n(x).
\end{align*}
 We next claim that there exist nonzero numbers $r_n$ and two polynomials $f(x)$ and $g(x)$ of degree at most two and one, respectively, not simultaneously zero, such that
$$
f_n(x)=r_n\, f(x), \quad g_n(x)=r_n\, g(x).
$$
Indeed, multiplying \eqref{3b} by $F_{n+1}$ and \eqref{3a} by $D_n(x)$,  we can eliminate $\mathcal{D}^2_q  P^{[k-1]}_{n}(x)$, and obtain, using \eqref{4} and shifting $n$ to $n-1$, the equation
\begin{align}\label{7}
f_{n-1}(x)\mathcal{D}^2_q  P^{[k-1]}_{n}(x)+g_{n-1}(x) \mathcal{S}_q \mathcal{D}_q  P^{[k-1]}_{n}(x)+\frac{\gamma_k}{\gamma_{n+1}}E_{n-1}P^{[k-1]}_{n}(x)=0.
\end{align}
(If $D_n(x)=0$, we combine directly \eqref{3b} and \eqref{4} to obtain \eqref{7}.)
Suppose that $E_{n-1}=0$, i.e.
$$
f_{n-1}(x)\mathcal{D}_q  P^{[k]}_{n-1}(x)=-g_{n-1}(x) \mathcal{S}_q P^{[k]}_{n-1}(x).
$$
Since $\mathcal{D}_q  P^{[k]}_{n-1}(x)$ and $\mathcal{S}_q P^{[k]}_{n-1}(x)$ have no common zeros, we have $f_{n-1}(x)=g_{n-1}(x)=0$, which is imposible according to \eqref{6} after shifting $n$ to $n-1$. Thus, by combining \eqref{6} and \eqref{7}, we can eliminate $P^{[k-1]}_{n}(x)$, and obtain the equation
\begin{align*}
&\left(f_n(x)-\frac{\gamma_{n+1}}{\gamma_{n-1}}\frac{\widetilde{E}_n}{E_{n-1}} C^{[k]}_{n-1} f_{n-1}(x)\right) \mathcal{D}_q  P^{[k]}_{n}(x)\\[7pt]
&\quad =-\left(g_n(x)-\frac{\gamma_{n+1}}{\gamma_{n-1}}\frac{\widetilde{E}_n}{E_{n-1}} C^{[k]}_{n-1} g_{n-1}(x)\right) \mathcal{S}_q P^{[k]}_{n}(x).
\end{align*}
Again, since $\mathcal{D}_q  P^{[k]}_{n}(x)$ and $\mathcal{S}_q P^{[k]}_{n}(x)$ have no common zeros, the desired conclusion follows. This allows us to rewrite \eqref{6} as
\begin{align}\label{final}
f(x)\mathcal{D}^2_q  P^{[k-1]}_{n}(x)+g(x) \mathcal{S}_q \mathcal{D}_q  P^{[k-1]}_{n}(x)+\lambda_n\,P^{[k-1]}_{n}(x)=0,
\end{align}
where $f(x)$ and $g(x)$ are polynomials of degree at most two and one, respectively, not simultaneously zero, and  $\lambda_n=\gamma_k/(r_n \gamma_{n-1})\widetilde{E}_n C^{[k]}_{n-1}\not=0$. Thus, by Theorem \ref{T}, $(P^{[k-1]}_{n})_{n\geq 0}$ are multiples of the Askey-Wilson polynomials, or special or limiting cases of them. Since now $(P_n)_{n\geq 0}$ and $(P^{[k-1]}_{n})_{n\geq 0}$  are two sequences of orthogonal polynomials, repeating the previous argument we conclude, for each $j=k-2, \dots, 0$, that $(P^{[j]}_{n})_{n\geq 0}$ are multiples of the Askey-Wilson polynomials, or special or limiting cases of them. The rest of the proof is trivial.

\section*{Acknowledgements}
The authors thank to Professor T. H. Koornwinder for helpful discussions and comments. This work was supported by the Centre for Mathematics of the University of Coimbra-UIDB/00324/2020, funded by the Portuguese Government through FCT/ MCTES. The second author thanks the support of the ERDF and Consejería de Economía, Conocimiento, Empresas y Universidad de la Junta de Andalucía (grant UAL18-FQM-B025-A).


\begin{thebibliography}{10}


\bibitem{KJ19} M. Kenfack Nangho and K. Jordaan, A characterization of Askey-Wilson polynomials, Proc. Amer. Math. Soc. 147 (2019) 2465-2480.




\bibitem{I05} {M. E. H. Ismail}, {Classical and quantum orthogonal polynomials in one variable. With two chapters by W. Van Assche.
With a foreword by R. Askey.}, Encyclopedia of Mathematics and its Applications \textbf{98}.
Cambridge University Press, Cambridge, 2005.

\end{thebibliography}
 \end{document}